\documentclass[11pt]{article}

\usepackage{amssymb,latexsym,amsfonts,verbatim,amscd}

\newtheorem{Def}{Definition}[section]
\newtheorem{Thm}[Def]{Theorem}
\newtheorem{Lem}[Def]{Lemma}

\newtheorem{Cor}[Def]{Corollary}

\newtheorem{Fac}[Def]{Fact}

\font\nat msbm10 scaled\magstephalf
\def\N{\hbox{\nat\char78}}

\def\telos{\hfill$\dashv$}

\font\goth eufm10

\begin{document}

\title{Some structural similarities between uncountable sets, powersets and the universe}
\author{Athanassios Tzouvaras}

\date{}
\maketitle

\begin{center}
Department  of Mathematics\\  Aristotle University of Thessaloniki \\
541 24 Thessaloniki, Greece \\
e-mail: \verb"tzouvara@math.auth.gr"
\end{center}

\begin{abstract}
We establish some  similarities/analogies between uncountable cardinals or powersets and the class $V$ of all sets. They concern mainly the Boolean algebras  ${\cal P}(\kappa)$, for a regular cardinal $\kappa$,  and ${\cal C}(V)$ (the class of subclasses of the universe $V$), endowed with some ideals, especially the ideal $[\kappa]^{<\kappa}$ for ${\cal P}(\kappa)$, and the ideal of sets $V$ for ${\cal C}(V)$.
\end{abstract}

\vskip 0.1in


\section{Introduction}
This article belongs roughly to the same line of thought as  \cite{Tz03} and \cite{Tz04b},  which on the one hand try to shed some light to the structure of powersets of infinite sets, and on the other, at philosophical level,   express  a skeptical attitude towards  the powerset axiom ($Pow$), that is the claim that the collection of subsets of an infinite set is  a set (rather than a proper class). To the same vein belongs also \cite{Tz10}, where it is proposed that the powerset and replacement axioms (especially the former one) would be better to be   replaced by  local (and thus relativized) versions of them that lead to a local version of the entire ${\rm ZFC}$. The approach of the present article differs, however,  in that we focus on the relationship and similarities between {\em uncountable} sets in general, rather than powersets alone, and the universe. Since the connection between the existence of powersets and the existence of uncountable sets seems  somewhat confused, let me elaborate briefly on this  point before  coming to the main body of the paper.

It is well-known that powersets and uncountable sets are closely related through the famous Cantor's theorem that every infinite powerset  is uncountable. However this connection seems to be only one-way: as is well-known,   the existence of uncountable sets is (relatively) consistent with the failure of $Pow$.  For example if $M\models{\rm ZFC}+2^{\aleph_0}=\aleph_2$,  then
$$H_{\omega_2}^M\models {\rm ZFC}^-+\mbox{``$\omega_1$ exists''}+\mbox{``${\cal P}(\omega)$ does not exist''},$$
where  $H_{\omega_2}=\{x:|TC(x)|<\aleph_2\}$ and ${\rm ZFC}^-={\rm ZFC}-\{Pow\}$. So, if $Unc$ denotes the existence of uncountable sets,  ${\rm ZFC}^-+Unc\not\vdash Pow$, i.e., $Pow$ is {\em formally} independent of $Unc$ over ${\rm ZFC}^-$.

The problem is that  a resident of $H_{\omega_2}^M$, in contrast to a resident of $M$,  has no clue as to {\em how} the uncountable $\omega_1$ has emerged. The resident of $M$ has a natural way to construct  $\omega_1$, as well as the other uncountable cardinals: they consider the collection $WO(\omega)$ of all well-orderings of $\omega$, which is a subcollection of ${\cal P}(\omega\times\omega)$. Since the latter is a {\em set} by $Pow$, so is $WO(\omega)$ by Separation. Then taking the set $WO(\omega)/\!\!\cong$ of isomorphism types of $WO(\omega)$, they arrive at  $\omega_1$. But for the resident of $H_{\omega_2}^M$, ${\cal P}(\omega\times\omega)$ is only a proper class. It is true that they can define $WO(\omega)$ in $H_{\omega_2}^M$, as a subclass of ${\cal P}(\omega\times\omega)$. But they cannot even define $WO(\omega)/\!\!\cong$, since each isomorphism type of a well-ordering of $\omega$ has continuum many members, so it is a proper class too.  Thus $\omega_1$ is  a completely {\em accidental} object of $H_{\omega_2}^M$, which seems to be there by mere chance.

If, alternatively, we start with  a model  $N\models{\rm ZFC}+2^{\aleph_0}=\aleph_1$, then
$$H_{\omega_2}^N\models {\rm ZFC}^-+\mbox{``${\cal P}(\omega)$ exists''}+\mbox{``${\cal P}({\cal P}(\omega))$ does not exist''},$$
and the problem now for a resident of $H_{\omega_2}^N$ is to  explain how the uncountable set ${\cal P}(\omega)$ has occurred, since the operation ${\cal P}(x)$ is not defined for all sets.

The role of $Pow$ in ${\rm ZF}$ is quite analogous to the role of exponentiation axiom ($Exp$) in (subsystems) of Peano arithmetic (${\rm PA}$). $Exp$ is the sentence (in the language of PA)  $\forall x,y\exists z(x^y=z)$, or more simply  $\forall x\exists y(2^x=y)$. If ${\rm I\Sigma}_n$ denotes the fragment of PA that consists of the elementary axioms of PA and the induction scheme restricted to $\Sigma_n$ formulas alone, then it is well-known that $Exp$ is provable in ${\rm I\Sigma}_1$ but unprovable in ${\rm I\Sigma}_0$. More precisely, ${\rm I\Sigma}_0<{\rm I\Sigma}_0+Exp<{\rm I\Sigma}_1$ (where for theories $T_1$, $T_2$, $T_1<T_2$ means that  $T_2$ is strictly stronger than $T_1$). Roughly, $Exp$ is with respect to ${\rm I\Sigma}_0$, what $Pow$ is with respect to ${\rm ZF}^-$, and  uncountable sets are analogous to ``very big'' numbers of arithmetic. The only difference is that   uncountable sets have a precise definition (sets of  size strictly larger than $\aleph_0$), while very big numbers do not.  Intuitively  $a$ is ``very big'' if  $a\geq 2^b$ for some  ``big'' $b$. As is well-known the models of ${\rm I\Sigma}_0$ are closed under $+$ and $\cdot$ but not, in general, under exponentiation. In particular, the operation $x^y$ is only partially defined in them.   Nevertheless, there are (nonstandard) models of ${\rm I\Sigma}_0+\neg Exp$ (see e.g. \cite{Pa71}, Theorems 4.3, 4.4) which, despite the failure of  $Exp$,  may contain  $2^a$ for some arbitrary  nonstandard number $a$. In the absence  of the full exponentiation operation, the  existence of $2^a$ can be seen by a resident of the model only as  {\em accidental}, exactly  as was seen the cardinal $\omega_1$ in the model $H_{\omega_2}^M$ of ${\rm ZFC}^-$ mentioned  above.

The conclusion of the preceding discussion is that the connection  between powersets and the uncountable is not just one-way. It is in fact a two-way dependence, if behind every  encounter with an uncountable set we seek  the reason and source  of its occurrence.

Having argued that  $Pow$ is, directly or indirectly, {\em responsible} for the existence of uncountable sets, a reasonable way to reinforce the skeptical attitude towards this axiom is to present concrete {\em mathematical facts} that establish  specific similarities/analogies  between  uncountable sets and the universe of  all sets $V$. This is what we are doing in this article. Namely,  we  compare the  Boolean algebra  ${\cal P}(\kappa)$ of an  infinite cardinal $\kappa$, with the Boolean algebra ${\cal C}(V)$ of subclasses of $V$, and point out certain similarities. The similarities  are of three kinds: (a) set-theoretic, (b) algebraic, and  (c) order-theoretic.  In the comparison of set-theoretic and algebraic   facts some ideals play a key role, namely the ideals $[\kappa]^{<\lambda}$, for $\lambda\leq \kappa$, and in particular  $[\kappa]^{<\kappa}$, for the case of cardinals, and the ideal of sets $V$ for ${\cal C}(V)$. We deal with each of the three kinds of similarities in  separate sections.

\vskip 0.1in

\textbf{Summary of contents.}
In section 2 we deal with a specific set-theoretic similarity. The main results are the following:

(a) The ideal $V$ of ${\cal C}(V)$ is {\em minimal} in a natural and  precisely defined way. Moreover, it is the {\em unique} minimal ideal of this kind (Theorem \ref{T:V-unique}).

(b) For every regular cardinal $\kappa$, the ideal $[\kappa]^{<\kappa}$ of  ${\cal P}(\kappa)$ is {\em minimal}  in the above sense, and also the {\em unique} minimal  ideal of this kind (Theorem \ref{T:regular}).

(c) (In sharp contrast to (b)), for every singular cardinal $\kappa$, there is {\em no} minimal ideal of ${\cal P}(\kappa)$ in the above sense (Theorem \ref{T:nominimal}).

In section 3 we deal with the algebraic similarities. Specifically  we compare the Boolean algebra  ${\cal P}(\kappa)$  equipped with an ideal $[\kappa]^{<\lambda}$, for any $\aleph_0\leq \lambda\leq \kappa$, with the Boolean algebra ${\cal C}(V)$ equipped with the ideal $V$.  Using Ehrenfeucht-Fra\"{i}ss\'{e} games, it is shown  that  these structures are elementarily equivalent, that is,  $\langle {\cal P}(\kappa),[\kappa]^{<\lambda}\rangle\equiv_{L_{b}(S)} \langle {\cal C}(V),V\rangle$, where $L_{b}(S)$ is the language of Boolean algebras augmented with a unary predicate $S(\cdot)$.

In section 4 we deal with the order-theoretic similarities. For convenience we identify $V$ with the class of ordinals $On$ by the help  of von Neumann's global choice $vN$, and compare the first-order structures $\langle \kappa,<\rangle$ and $\langle On,<\rangle$, where $<$ is the natural well-orderings of $On$ and $\kappa$ (i.e., $\in$). Extending  results of  \cite{DMT78} and \cite{Tz04a} we show  that $\langle \kappa,<\rangle\equiv_{L_{ord}} \langle On,<\rangle$, where $L_{ord}$ is the first-order language of ordering. This equivalence is further strengthened to $\langle \kappa,{<,}{\rm Cof}(\kappa)\rangle\equiv^{\forall_1^1}_{pos} \langle On,<,{\rm Cof}(On)\rangle$, where ${\rm Cof}(\kappa)$ (resp. ${\rm Cof}(On)$) denotes the class of cofinal subsets of $\langle \kappa,<\rangle$ (resp. cofinal subclasses of $\langle On,<\rangle$) and $\equiv^{\forall_1^1}_{pos}$ denotes elementary equivalence with respect  to all positive $\forall_1^1$ formulas of the monadic (second-order) language $L_{mon}$, which is $L_{ord}$ augmented with set variables and $\in$.

\vskip 0.1in

\textbf{Some  notational conventions.}
In order to talk about subclasses of $V$ we  need to work in a class theory (rather than  ${\rm ZFC}$), like G\"{o}del-Bernays theory with choice ${\rm GBC}$. So throughout our discussion  will be held having as ground theory either ${\rm GBC}$ or ${\rm ZFC}$. Concerning the collection  ${\cal C}(V)$  of all subclasses of $V$ when working in ${\rm GBC}$, this  is actually a ``virtual'', i.e.  syntactic, entity, since in order to accommodate it literally as an extensional object, we would need a third-order set theory. However this is not a serious problem.  We can always treat ${\cal C}(V)$ in model-theoretic terms, that is, think that we work in an arbitrary model $\mbox{\goth M}=\langle M,{\cal M}\rangle$ of ${\rm GBC}$, where  $M$ contains  the sets and ${\cal M}$ contains the classes of $\mbox{\goth M}$. Then $V^{\mbox{\goth M}}=M$ while ${\cal C}(V)^{\mbox{\goth M}}={\cal M}$. Thus ${\cal C}(V)$ is simply the range of  class variables and,  like ${\cal P}(A)$, an atomic Boolean algebra. On the other hand,  when we talk about  {\em ideals} of ${\cal C}(V)$, which are  collections ${\cal I}\subseteq {\cal C}(V)$, these will always be {\em definable}  by a formula $\phi(X)$ with one free class variable.

We have already used above the notation $[\kappa]^{<\lambda}$. In general, for any  set $A$ and cardinal $\kappa\leq |A|$ we let
$$[A]^{<\kappa}=\{x\subseteq A:|x|<\kappa\}\footnote
{This set is usually denoted also ${\cal P}_\kappa(A)$.}, \ \
[A]^{\leq \lambda}=\{x\subseteq A:|x|\leq \lambda\}.$$
Also, in view of  $AC$,  we can identify every powerset ${\cal P}(A)$ with a cardinal $2^\kappa$, and write ${\cal P}(2^\kappa)$ instead of ${\cal P}({\cal P}(A))$, $[2^\kappa]^{<\lambda}$ instead of $[{\cal P}(A)]^{<\lambda}$, and so on.

\section{A set-theoretic similarity between  $V$ and regular cardinals: minimality and uniqueness of access ideals}
Throughout  upper case letters $X$, $Y$ are used as class variables for the language of ${\rm GBC}$,  while lower case letters $x$, $y$ denote sets (elements of classes). For two  classes $X,Y$, let $|X|=|Y|$ denote the fact that $X$ and $Y$ are equipollent, i.e. there is a bijection between them. The axiom  $vN$ says that all proper classes are equipollent, or equivalently $|V|=|On|$. Following the conventions agreed in the Introduction we shall write $2^{\aleph_0}$ instead of ${\cal P}(\omega)$ and in general $2^\kappa$ instead of ${\cal P}(\kappa)$.

Definable classes are produced as usual by formulas $\phi(x)$ with a set variable, i.e., $X=\{x:\phi(x)\}$. But also  classes of classes are  permitted, provided they are defined by formulas $\phi(X)$ with a class free variable. We  denote such classes by calligraphic letters  ${\cal I}$, ${\cal J}$ and treat them as syntactic objects, namely ${\cal I}$, ${\cal J}$ can be thought of as formulas, and writing  ${\cal I}\subseteq {\cal J}$ we just mean that $(\forall X)({\cal I}(X)\rightarrow {\cal J}(X))$. For any such ${\cal I}$ we have ${\cal I}\subseteq {\cal C}(V)$, where ${\cal C}(V)$ is the class of all subclasses of $V$.
Given a class $X$, an {\em ideal} of ${\cal C}(X)$ (or an ideal on $X$) is a  class ${\cal I}\subseteq {\cal C}(X)$ such that:

(i) If $Y_1,Y_2\in {\cal I}$, then $Y_1\cup Y_2\in {\cal I}$,

(ii) If $Y\in {\cal I}$ and $Z\subseteq Y$, then $Z\in {\cal I}$,

(iii) ${\cal I}$ is proper, i.e., $X\notin {\cal I}$.

\noindent Trivially  $V$ is an ideal of ${\cal C}(V)$.

\subsection{Accessibility of classes and sets by ideals}
A (total) {\em preordering} on a class $X$ is a binary relation $\preccurlyeq$ which is  reflexive, transitive and total, where the latter means that for all $x,y\in X$ $x\preccurlyeq y$ or $y\preccurlyeq x$. $x\prec y$ means $x\preccurlyeq y$ and $x\neq y$.  If we add the antisymmetric property, the preordering becomes a linear ordering. To avoid trivialities we shall assume throughout that no preordered  class $\langle X,\preccurlyeq\rangle$ has a last element. Given a preordering $\langle X,\preccurlyeq\rangle$  and  $x\in X$,  ${\preccurlyeq_x}$ denotes the initial segment of $X$ determined by $x$, i.e., ${\preccurlyeq_x}=\{y\in X: y\preccurlyeq x\}$. A  $Y\subseteq X$ is said to be {\em bounded} w.r.t.~$\preccurlyeq$ if there is $x\in X$ such that $Y\subseteq {\preccurlyeq_x}$. An unbounded $Y\subseteq X$ is called also {\em cofinal}.  For every preordering $\preccurlyeq$ of  $X$, the bounded subsets of $X$ form a natural ideal of  ${\cal C}(X)$ denoted $Seg(\preccurlyeq)$, that is
$$Seg(\preccurlyeq)=\{Y\subset X:\exists x\in X(Y\subseteq{\preccurlyeq_x})\}.$$
The fact that $\langle X,\preccurlyeq\rangle$ has no last element guarantees that the ideal  $Seg(\preccurlyeq)$ is proper.

\begin{Def} \label{D:access}
{\em (GBC or ZFC) Let ${\cal I}$ be an ideal on a class or set $X$.

(1) ${\cal I}$ is said to be an} access ideal for $X$,  {\em or} $X$ is ${\cal I}$-accessible, {\em  if there is a preordering  $\preccurlyeq$ of $X$ such that   $Seg(\preccurlyeq)\subseteq {\cal I}$. In that case  $\preccurlyeq$ is said to be an} ${\cal I}$-preordering.

{\em (2) The ideal ${\cal I}$ is said to be a} minimal access ideal, {\em or just} minimal for $X$,  {\em if  it is an access ideal for $X$ and for any ideal ${\cal J}\subsetneq {\cal I}$, $X$ is not ${\cal J}$-accessible. }
\end{Def}

\noindent Trivially every preordering  $\preccurlyeq$ of $X$ is a $Seg(\preccurlyeq)$-preordering. Notice that, although the preceding definitions are expressible in  ${\rm GBC}$ or ${\rm ZFC}$ (when $X$ and ${\cal  I}$ are definable objects), it is not to be meant that it is also {\em decidable} in these theories whether, for given   $X$ and ${\cal I}$, (a) $X$ is ${\cal I}$-accessible or not, and (b) ${\cal I}$ is minimal for $X$ or not. These are relative facts, in general, depending on the model of ${\rm GBC}$/${\rm ZFC}$ we work in. For example, as we shall see below, the cardinal $2^{\aleph_0}$ is accessible by  the ideal $[2^{\aleph_0}]^{\leq \aleph_0}$ of its countable subsets if and only if $2^{\aleph_0}=\aleph_1$. Moreover, in the latter case $[2^{\aleph_0}]^{\leq \aleph_0}$ is the unique minimal access ideal for $2^{\aleph_0}$. More generally, if $2^{\aleph_0}=\aleph_\alpha$ holds in a model $\mbox{\goth M}$ and $\aleph_\alpha$ is regular, then $[2^{\aleph_0}]^{< \aleph_\alpha}$ is the unique minimal access ideal for $2^{\aleph_0}$ in $\mbox{\goth M}$. In contrast, if $\aleph_\alpha$ is singular, then $2^{\aleph_0}$ is $[2^{\aleph_0}]^{< \aleph_\alpha}$-accessible but there is no minimal access ideal for  $2^{\aleph_0}$.

On the other hand,  the class $V$ is $V$-accessible {\em independently} of the underlying model, in view of the axiom of regularity and the induced $V$-preordering $x\preccurlyeq_r y\Leftrightarrow rank(x)\leq rank(y)$ of the universe.\footnote{Actually this fact is the  reason  that in the definition of accessibility \ref{D:access} given above we use  preorderings rather than total orderings or well-orderings.} Similarly, for every {\em regular} cardinal $\kappa$, $\kappa$ is provably $[\kappa]^{<\kappa}$-accessible. More important are the facts shown below, that the ideals $V$ and $[\kappa]^{<\kappa}$ are both {\em provably} (a) {\em minimal} access ideals,  and (b) the {\em unique} minimal access ideals for $V$ and any regular $\kappa$, respectively.

Needles to say that for  $V$, as a well as for any cardinal   $\kappa$, there is an abundance of access ideals that arise  in the form of $Seg(\preccurlyeq)$ for the various preorderings $\preccurlyeq$ of $V$ or $\kappa$,  so the proof of existence of  {\em minimal} such ideals and, further, of {\em unique minimal} ones  is by no means trivial.  For example  partition $\kappa$ into  finitely  or infinitely many  pieces and then rearrange the pieces into an arbitrary  total ordering in which each piece keeps its  standard ordering.  Then  a total ordering of $\kappa$  arises and hence an access ideal for it.

According to Definition \ref{D:access}, if a class $X$ is ${\cal I}$-accessible for an  ideal ${\cal I}$ within a model $\mbox{\goth M}$, then there is a  preordering  $\preccurlyeq$ of $X$ in $\mbox{\goth M}$ such that every  $\preccurlyeq$-bounded set  of $X$ belongs to ${\cal I}$, i.e.,
\begin{equation} \label{E:oneway}
\mbox{\goth M}\models (\forall Y\subseteq X)(Y \  \mbox{is $\preccurlyeq$-bounded} \rightarrow Y\in {\cal I}).
\end{equation}
This definition obviously allows the ideal ${\cal I}$ to contain also unbounded subsets of $X$. If, however, ${\cal I}$ is minimal, then the converse of (\ref{E:oneway}) holds too.  Namely the following holds.

\begin{Lem} \label{L:minimal}
Let  $\mbox{\goth M}\models{\rm GBC}$ and let ${\cal I}$ be a minimal access ideal for a  definable class $X$ in $\mbox{\goth M}$.  Then for every ${\cal I}$-preordering $\preccurlyeq\in \mbox{\goth M}$ of $X$:
$$\mbox{\goth M}\models (\forall Y\subseteq X)(Y \ \mbox{is $\preccurlyeq$-bounded} \ \leftrightarrow \ Y\in {\cal I}).$$
\end{Lem}

{\em Proof.} One direction follows from Definition \ref{D:access}.  For the converse, assume that there is an ${\cal I}$-preordering $\preccurlyeq\in \mbox{\goth M}$ of $X$ such that ${\cal I}$ contains a set  $Y_0\subseteq X$  unbounded in $\langle X,\preccurlyeq\rangle$.  Then  ${\cal I}$ is not a minimal access ideal for $X$.  Because   $X$ is $Seg(\preccurlyeq)$-accessible, and $Seg(\preccurlyeq)\subsetneq {\cal I}$, since  $Y_0\in {\cal I}\backslash Seg(\preccurlyeq)$.  \telos

\vskip 0.2in

\underline{\textbf{Minimality and uniqueness of $V$ as an access ideal for $V$}}

\begin{Thm} \label{T:V-minimal}
{\rm (GBC)} The class  $V$ is a minimal access ideal for  $V$.
\end{Thm}

{\em Proof.} By the regularity axiom of ${\rm GBC}$, the relation  $x\preccurlyeq y$ iff  $rank(x)\leq rank(y)$ is a $V$-preordering, thus $V$ is $V$-accessible in ${\rm GBC}$. To prove the minimality of $V$, let ${\cal J}$ be any ideal on $V$ such that ${\cal J}\subsetneq V$. It suffices to see that $V$ is not ${\cal J}$-accessible. Assume the contrary, and let $\preccurlyeq_1$ be a  preordering of $V$ such that for every $x\in V$, $\preccurlyeq_{1,x}\in J$. Pick a $z\in V\backslash {\cal J}$. If  $z$ were cofinal in $\langle V,\preccurlyeq_1\rangle$, then  $V=\bigcup_{x\in z}\preccurlyeq_{1,x}$, and since for every $x$, $\preccurlyeq_{1,x}\in {\cal J}\subseteq V$, by Replacement and Union $\bigcup_{x\in z}\preccurlyeq_{1,x}=V$ would be a  set, a contradiction.  If, on the other hand, $z$ were not cofinal in $\langle V,\preccurlyeq_1\rangle$, we should  have $z\subseteq \preccurlyeq_{1,x}$ for some $x$, so $z\in {\cal J}$ since $ \preccurlyeq_{1,x}\in {\cal J}$, a contradiction again.  \telos

\vskip 0.2in

In fact $V$ is not just a minimal access ideal for $V$, but it is the {\em unique} minimal access ideal. We shall need, first, the following simple consequence of the axiom of choice $AC$. This is rather folklore but we include its proof here for the reader's convenience.

\begin{Lem} \label{L:around}
{\rm (GBC)} Let $X$ be a proper class. Then for every cardinal $\kappa$ there is set $b\subset X$ such that $|b|=\kappa$.
\end{Lem}

{\em Proof.} Given the proper class $X$, let $R=\{rank(x):x\in X\}$.  $R$ is a proper class of ordinals and for every such class there is  in GB a strictly increasing  enumeration  $R=\{\alpha_\xi:\xi\in On\}$, i.e., a bijection  $F:On\rightarrow R$ such that $\xi<\zeta\Leftrightarrow F(\xi)<F(z)$. (This can be shown as follows: let $S$ be the class of set functions $f$ such that $dom(f)=\xi\in On$, $rng(f)\subset R$ and for every $\delta\in dom(f)$, $f(\delta)=\min(R\backslash \{f(\gamma):\gamma<\delta\})$. $S$ is defined by a normal formula so it is a class of GB. Then $F:On\rightarrow R$ is defined as follows: $F(\xi)=\alpha\Leftrightarrow (\exists f\in S)(\xi\in dom(f) \ \wedge \ f(\xi)=\alpha)$.) If $d_\xi=X\cap (V_{\alpha_\xi+1}\backslash V_{\alpha_\xi})$, then $\{d_\xi:\xi\in On\}$, is a proper class of nonempty disjoint sets. Let  $s=\{d_\xi:\xi<\kappa\}$. $s$ is a set of nonempty sets, so by $AC$ it has a choice function $f$. If $b=f[s]$, then $b$ is as required. \telos

\begin{Thm} \label{T:V-unique}
{\rm (GBC)} $V$ is the unique  minimal access ideal for $V$.
\end{Thm}

{\em Proof.} Towards reaching a contradiction assume that  ${\cal I}\neq V$ is another  minimal access ideal for  $V$. Then we cannot have $V\subseteq {\cal I}$, for since ${\cal I}\neq V$, that would mean that  $V\subsetneq{\cal I}$, which contradicts the minimality of ${\cal I}$. Thus  $V\not\subseteq {\cal I}$. Pick  a set $a\in V\backslash {\cal I}$ and fix an ${\cal I}$-preordering  $\preccurlyeq$ of $V$. (Remember that all preorderings used throughout  have no last element.) It suffices to  construct a preordering $\preccurlyeq_1$ of $V$ such that $Seg(\preccurlyeq_1)\subsetneq{\cal I}$, contradicting thus  the minimality of ${\cal I}$.

Notice first that the set $a\in V\backslash {\cal I}$ must be cofinal in $\langle V,\preccurlyeq\rangle$, otherwise $a\subseteq {\preccurlyeq_x}$ for some $x$, and hence $a\in {\cal I}$ since ${\preccurlyeq_x}\in {\cal I}$.  Also for some $x\in a$, ${\preccurlyeq_x}$ must be  a proper class,  otherwise $V=\bigcup_{x\in a}{\preccurlyeq_x}$ would be a set. Fix  an element $x_0\in a$ such  that $\preccurlyeq_{x_0}$ is a proper class and let $X=\preccurlyeq_{x_0}$. Let also $Y=V\backslash X$ be the final segment of $(V,\preccurlyeq)$ determined by $x_0$. Then we can write
\begin{equation} \label{E:V-analysis1}
V=X\uplus Y
\end{equation}
for the ordered disjoint union of the initial and final segments to which $x_0$  splits $\langle V,\preccurlyeq\rangle$. We shall construct a preordering $\preccurlyeq_1$ such that if $J=Seg(\preccurlyeq_1)$, then $J\subsetneq {\cal I}$.  Let us set $a^*=a\cap Y$. Since $a$ is cofinal in $\preccurlyeq$ and infinite (because $\preccurlyeq$ does not have a last element),  so is $a^*$.  Let $|a^*|=\kappa$ and let  ${a^*}=\{y_\xi:\xi<\kappa\}$. Using Lemma \ref{L:around}, pick a subset $b$ of the proper class $X$ such that $|b|=|a^*|=\kappa$ and let $b=\{z_\xi:\xi<\kappa\}$. The idea is to cut $b$ off $X$ and make it cofinal in the new preordering $\preccurlyeq_1$, retaining $\preccurlyeq$ to the rest parts of $V$. This can be  obtained if we define $\preccurlyeq_1$ so that $\preccurlyeq_1$ agrees with $\preccurlyeq$ on $V\backslash b$, while $b$ and  $a^*$  both run cofinally to $V$ (and to each other) with respect to $\preccurlyeq_1$. So let
$$R=[\preccurlyeq\cap (V\backslash b)^2]\cup\{\langle y_\xi,z_\xi\rangle, \langle z_\xi,y_\xi\rangle:\xi<\kappa\},$$
and let $\preccurlyeq_1$ be the transitive closure of the binary relation  $R$, i.e., $\langle x,z\rangle\in \preccurlyeq_1$ whenever there is $y$ such that $\langle x,y\rangle$ and $\langle y,z\rangle$ belong to $\preccurlyeq_1$. By this definition $y_\xi\preccurlyeq_1z_\xi$ and $z_\xi\preccurlyeq_1y_\xi$ for all $\xi<\kappa$, and  $\preccurlyeq_1$ is a preordering of $V$ which,  in analogy  to  (\ref{E:V-analysis1}),  splits  $V$ into  the ordered union
\begin{equation} \label{E:V-analysis2}
V=(X\backslash b)\uplus (Y\cup b),
\end{equation}
with   $X\backslash b$, $Y$ and  $b$ being pairwise disjoint. Moreover $a^*\subseteq Y$, hence $a^*\cap b=\emptyset$ and $a^*$ is cofinal in $\langle V,\preccurlyeq_1\rangle$. Let
$${\cal J}=Seg(\preccurlyeq_1)=\{Z:(\exists w)(Z\subseteq \preccurlyeq_{1,w})\}.$$
We show first that ${\cal J}\subseteq {\cal I}$.   Let $Z\in {\cal J}$ and let  $Z\subseteq \preccurlyeq_{1,w}$ for some $w$. We shall show that $Z\in {\cal I}$.  By the picture of $V$ in (\ref{E:V-analysis2}) as an ordered union, $Z$ is written  as a disjoint union $Z=X_1\cup X_2\cup x_3$, where $X_1\subseteq X\backslash b$, $X_2\subseteq Y$ and $x_3\subseteq b$. Since $X\backslash b$ and $b$ already belong to ${\cal I}$, it follows that $X_1,x_3\in {\cal I}$. It remains to show that $X_2\in {\cal I}$. We have $X_2\subseteq Y$ and  $X_2\subseteq \preccurlyeq_{1,w}$. Clearly either  $w\in Y$ or $w\in b$. If  $w\in Y$, then since $\preccurlyeq$ and $\preccurlyeq_1$ agree on $Y$, we have  $X_2\subseteq \preccurlyeq_w$, so $X_2\in {\cal I}$ because  $\preccurlyeq$ is an ${\cal I}$-preordering.  Assume  $w\in b$, so $w=z_\xi$ for some $\xi<\kappa$. By definition $z_\xi\preccurlyeq_1y_\xi$, where  $y_\xi\in a^*\subseteq Y$, i.e.,  $y_\xi\in Y$. Then   $X_2\subseteq \preccurlyeq_{1,z_\xi}\subseteq \preccurlyeq_{1,y_\xi}$. Since $y_\xi\in Y$ as before  $X_2\subseteq \preccurlyeq_{y_\xi}$, whence again  $X_2\in {\cal I}$. Thus all three pieces $X_1,X_2,x_3$ of $Z$ belong to ${\cal I}$, and hence  $Z\in {\cal I}$. This proves  ${\cal J}\subseteq {\cal I}$. On the other hand $b\notin {\cal J}$ because it is cofinal in $\preccurlyeq_1$, while $b\in {\cal I}$ as a subset of $X$ which belongs to ${\cal I}$. Therefore ${\cal J}\subsetneq {\cal I}$. This contradicts the minimality of  ${\cal I}$ and completes the proof. \telos

\vskip 0.2in

\underline{\textbf{Minimality and uniqueness of access ideals on cardinals}}

Before considering access ideals on powersets, i.e., on cardinals $2^\kappa$, we must consider access  ideals on general uncountable cardinals. We shall see that there is a sharp difference between regular and singular cardinals with respect to this matter. We deal first with regular cardinals.

\begin{Thm} \label{T:regular}
{\rm (ZFC)} Let $\kappa$ be an uncountable cardinal and $<$ be its natural well-ordering. Then:

(a) $Seg(<)=[\kappa]^{<\kappa}$  iff $\kappa$ is regular.

(b) If $\kappa$ is regular, then  $[\kappa]^{<\kappa}$ is a minimal access ideal for $\kappa$.

(c) Moreover if $\kappa$ is regular, $[\kappa]^{<\kappa}$ is the unique  minimal access  ideal for $\kappa$.
\end{Thm}

{\em Proof.} (a) Firstly,  for every $\kappa$, $Seg(<)\subseteq [\kappa]^{<\kappa}$.  Suppose $\kappa$ is regular and $x\in [\kappa]^{<\kappa}$. Then $x$ is bounded in $\langle \kappa,<\rangle$, thus $x\in  Seg(<)$, so $[\kappa]^{<\kappa}\subseteq Seg(<)$ and therefore $Seg(<)=[\kappa]^{<\kappa}$. Conversely, if $\kappa$ is not regular there is $x\subset \kappa$ cofinal in $\langle \kappa,<\rangle$ with $|x|<\kappa$. Thus $x\in [\kappa]^{<\kappa}\backslash Seg(<)$ and hence $Seg(<)\subsetneq [\kappa]^{<\kappa}$.

(b) Let $\kappa$ be regular and assume $[\kappa]^{<\kappa}$ is not minimal for $\kappa$, i.e., there is an ideal ${\cal I}\subsetneq [\kappa]^{<\kappa}$ and an ${\cal I}$-preordering $\preccurlyeq$ for $\kappa$. Pick $x\in [\kappa]^{<\kappa}\backslash {\cal I}$. Then $x$ is cofinal in $\langle \kappa, \preccurlyeq\rangle$, because otherwise $x\subseteq \preccurlyeq_\xi\in {\cal I}$, for some $\xi\in\kappa$, and hence $x\in {\cal I}$, contradiction. This implies that $\kappa=\bigcup_{\xi\in x}\preccurlyeq_\xi$. But for each $\xi$, $|\preccurlyeq_\xi|<\kappa$, since $\preccurlyeq_\xi\in {\cal I}\subseteq [\kappa]^{<\kappa}$, and also $|x|<\kappa$. This clearly contradicts the regularity for $\kappa$.

(c) Let $\kappa$ be regular. The proof that $[\kappa]^{<\kappa}$  is the unique  minimal ideal for $\kappa$ is quite similar to  the proof of Theorem  \ref{T:V-unique} for the corresponding fact about $V$, so we simply sketch it. Let ${\cal I}\neq [\kappa]^{<\kappa}$ be another  minimal access ideal for $\kappa$ and fix an ${\cal I}$-preordering  $\preccurlyeq$ of $\kappa$. If $[\kappa]^{<\kappa}\subseteq {\cal I}$, then $[\kappa]^{<\kappa}\subsetneq{\cal I}$, so ${\cal I}$ would not be minimal. Therefore   $[\kappa]^{<\kappa}\not\subseteq {\cal I}$.  Pick $a\in [\kappa]^{<\kappa}\backslash {\cal I}$. Then, as in the proof of \ref{T:V-unique},  $a$ is cofinal in $\langle\kappa,\preccurlyeq\rangle$, and for some $x\in a$, $|{\preccurlyeq_x}|=\kappa$. Because $\kappa=\bigcup_{x\in a} {\preccurlyeq_x}$, and $|a|<\kappa$, so if $|{\preccurlyeq_x}|<\kappa$ for all $x\in\kappa$, $\kappa$ would be singular. Fix such a $x_0\in \kappa$ and let $X=\preccurlyeq_{x_0}$, so $|X|=\kappa$,  $Y=\kappa\backslash X$, and $a^*=a\cap Y$. Let $|a^*|=\lambda<\kappa$, and pick $b\subset X$, with $|b|=|a^*|=\lambda$. Let also $a^*=\{y_\xi:\xi<\lambda\}$ and $|b|=\{z_\xi:\xi<\lambda\}$. We define a preordering $\preccurlyeq_1$ of $\kappa$ in complete analogy to that of \ref{T:V-unique}, making $b$ cofinal in $\langle\kappa,\preccurlyeq_1\rangle$, and so that $\kappa$ is written  $$\kappa=(X\backslash b)\uplus (Y\cup b),$$
with   $X\backslash b$, $Y$ and  $b$ being pairwise disjoint. If ${\cal J}=Seg(\preccurlyeq_1)$, it is shown as before that ${\cal J}\subsetneq {\cal I}$, which contradicts the minimality of ${\cal I}$.  \telos

\vskip 0.2in

Now when coming to powersets, Theorem \ref{T:regular} (c) gives its place to the following variant that claims  consistency rather than provability.

\begin{Thm} \label{T:uni-minimal}
{\rm (ZFC)} For any $\kappa\geq \aleph_0$ and any regular $\lambda>\kappa$, it is (relatively) consistent to believe that $[2^\kappa]^{<\lambda}$ is the unique minimal access ideal on $2^\kappa$.
\end{Thm}

{\em Proof.} It is well-known that the following is consistent with ${\rm ZFC}$: $\lambda>\kappa$ is regular and  $2^\kappa=\lambda$. This is obtained by a standard Cohen forcing extension $\mbox{\goth M}$ in which $\lambda=2^\kappa$ is still regular, so  it follows from \ref{T:regular} (c) that in $\mbox{\goth M}$, $[2^\kappa]^{<\lambda}=[\lambda]^{<\lambda}$ is the unique minimal access ideal on $2^\kappa$. \telos

\vskip 0.2in

We come next to the ${\cal I}$-accessibility of singular cardinals. First recall that for every preordered set  $\langle A,\preccurlyeq\rangle$, the {\em cofinality} of  $\langle A,\preccurlyeq\rangle$, denoted ${\rm cof}(A,\preccurlyeq)$,  is the least cardinal $\lambda$ such that: There is a set $X\subseteq A$ such that $|X|=\lambda$, and $\langle X,\preccurlyeq\rangle$ is well-ordered and cofinal to $\langle A,\preccurlyeq\rangle$. Clearly, for every infinite set preordering $\langle A,\preccurlyeq\rangle$ having no last element, ${\rm cof}(A,\preccurlyeq)$ is an infinite regular cardinal. Therefore, as an immediate consequence we have the following.

\begin{Fac} \label{F:nocof}
If  $\kappa$ is a singular cardinal, then for every preordering  $\preccurlyeq$ of $\kappa$, there is a regular $\lambda<\kappa$ such that ${\rm cof}(\kappa,\preccurlyeq)=\lambda$.
\end{Fac}

Note  that for any $\kappa$, either regular or singular, ${\rm cof}(\kappa,\preccurlyeq)$ need not be equal to  ${\rm cf}(\kappa)$. ${\rm cof}(\kappa,\preccurlyeq)$ depends heavily on  $\preccurlyeq$ and may be less  or greater  than ${\rm cf}(\kappa)$. For example for every $\kappa$ there is $\preccurlyeq$ such that ${\rm cof}(\kappa,\preccurlyeq)=\aleph_0$.  On the other hand  ${\rm cf}(\aleph_{\omega_1})=\aleph_1$, but we can easily find  orderings of $\aleph_{\omega_1}$ with cofinalities $\aleph_n>\aleph_1$. More generally, given any regular $\lambda$ and $\kappa>\lambda$, we can find a  preordering $\preccurlyeq_\lambda$ of $\kappa$ such that ${\rm cof}(\kappa,\preccurlyeq_\lambda)=\lambda$. (This can be done by splitting $\kappa$ into $\lambda$ (disjoint) pieces $x_i$, $i<\lambda$, e.g. by the help of a surjection $f:\kappa\rightarrow \lambda$,  and then define $\preccurlyeq_\lambda$ by totally ordering the pieces $x_i$ so that $x_i \prec_\lambda x_j$ iff $i<j$.)

Next we show  that, even  when $\kappa$ is singular, for every $\lambda<\kappa$, $\kappa$ is not $[\kappa]^{<\lambda}$-accessible, i.e., there is no
$[\kappa]^{<\lambda}$-preordering of $\kappa$. This is due to the following fact.

\begin{Lem} \label{L:anysegment}
Let $\kappa$ be an uncountable cardinal and let $\aleph_0\leq \lambda<\kappa$. For every preordering $\preccurlyeq$ of $\kappa$, there is an $x\in\kappa$ such that $|{\preccurlyeq_x}|\geq \lambda$.
\end{Lem}

{\em Proof.} Assume for contradiction that for a given $\lambda<\kappa$, there is a preordering $\preccurlyeq$ of $\kappa$ such that  $|{\preccurlyeq_x}|<\lambda$ for every $x\in \kappa$.  For each $x\in \kappa$, let $[x]=\{y\in\kappa:x\preccurlyeq y \ \& \ y\preccurlyeq x\}$ and let $\kappa/\!\!\sim=\{[x]:x\in \kappa\}$.  Then $\preccurlyeq$ is (essentially) a  linear ordering of $\kappa/\!\!\sim$. Since $|{\preccurlyeq_x}|< \lambda$ for all $x$, we have that $|[x]|<\lambda$ for all $x\in\kappa$, and therefore $|\kappa/\!\!\sim|=\kappa$. So it suffices to prove the claim for the linearly  ordered set $\langle\kappa/\!\!\sim,\preccurlyeq\rangle$. For simplicity we pick an element $s_x$  from each class $[x]$, and we set  $S=\{s_x:x\in \kappa\}$. Then  $|S|=\kappa$, and we prove the claim for the isomorphic set  $\langle S,\preccurlyeq\rangle$. Now either $S$ contains  a cofinal subset $Y$ of cardinality $\lambda$, or not. In the first case  $S=\bigcup_{s\in Y} \preccurlyeq_s$, where $|Y|=\lambda$ and $|\preccurlyeq_s|<\lambda$, for every $s\in Y$. But then $|S|=\lambda$, a contradiction. In the other case,  $S$ does not contain any cofinal set of cardinality $\lambda$. So  if we pick a $Y\subseteq S$ with $|Y|=\lambda$, then there is an $s\in S$ such that $Y\subseteq \preccurlyeq_s$. But then $|\preccurlyeq_s|\geq \lambda$, a contradiction again. \telos

\begin{Cor} \label{C:noreach}
Let $\aleph_0\leq \lambda<\kappa$. Then for every ${\cal I}\subseteq [\kappa]^{<\lambda}$, $\kappa$ is not ${\cal I}$-accessible, i.e., there is no  ${\cal I}$-preordering of $\kappa$.
\end{Cor}

{\em Proof.} If $\preccurlyeq$ were a $[\kappa]^{<\lambda}$-preordering of $\kappa$, then $|{\preccurlyeq_x}|<\lambda$ for every $x\in\kappa$, which contradicts \ref{L:anysegment}. \telos

\vskip 0.2in

Next theorem is in sharp contrast to Theorem \ref{T:regular} (c).

\begin{Thm} \label{T:nominimal}
{\rm (ZFC)} If  $\kappa$ is a singular cardinal, then there is no minimal access ideal on $\kappa$.  In particular, $Seg(<)$ is not minimal for $\kappa$.
\end{Thm}

{\em Proof.} Let $\kappa$ be singular. We have to show that  for every preordering $\preccurlyeq$ of $\kappa$, $Seg(\preccurlyeq)$ is not minimal. That is, there is a preordering $\preccurlyeq_1$ of $\kappa$ such that $Seg(\preccurlyeq_1)\subsetneq Seg(\preccurlyeq)$. It is remarkable that the construction of $\preccurlyeq_1$ follows exactly the same steps as the corresponding constructions in the proofs  of  theorems \ref{T:V-unique} and  \ref{T:regular} (c), although the latter theorems make claims to the opposite direction from  that of the present one.

Fix a preordering  $\preccurlyeq$  of  $\kappa$ and let $\lambda={\rm cof}(\kappa,\preccurlyeq)$. By Fact \ref{F:nocof}, $\lambda<\kappa$. Also by Lemma \ref{L:anysegment} there is $x_0\in \kappa$ with $|\preccurlyeq_{x_0}|\geq \lambda$.
We construct $\preccurlyeq_1$ imitating the proofs of \ref{T:V-unique} and \ref{T:regular} (c). Namely, let $X=\preccurlyeq_{x_0}$, let $b\subseteq X$ with $|b|=\lambda$, and let $Y=\kappa\backslash X$. Pick also a set $a$ with $|a|=\lambda$,  which is well-ordered by $\preccurlyeq$ with order-type $\lambda$, and $\langle a,\preccurlyeq\rangle$ is cofinal to $\langle \kappa,\preccurlyeq\rangle$. So we can write $a=\{\alpha_\xi:\xi<\lambda\}$, where $\alpha_\xi\prec \alpha_\zeta \Leftrightarrow \xi<\zeta$. Without loss of generality we may assume that $a\subseteq Y$. Let $b=\{\beta_\xi:\xi<\lambda\}$. As in the previous proofs we cut off $b$ from $X$ and make it cofinal to $a$. Namely, we let
 $$R=[\preccurlyeq \cap (\kappa\backslash b)^2]\cup\{\langle \alpha_\xi,\beta_\xi\rangle, \langle \beta_\xi,\alpha_\xi\rangle:\xi<\lambda\}. $$
Also we  set  $\preccurlyeq_1=$ the transitive closure of the binary relation $R$. Then $\preccurlyeq_1$ is  as required. For  exactly as in the proof of \ref{T:V-unique} we see  that $Seg(\preccurlyeq_1)\subseteq Seg(\preccurlyeq)$, while on the other hand $b\in Seg(\preccurlyeq)\backslash Seg(\preccurlyeq_1)$. Thus $Seg(\preccurlyeq_1)\subsetneq Seg(\preccurlyeq)$ and this completes the proof. \telos

\vskip 0.2in

As a complement to  Theorem \ref{T:uni-minimal} we have the following.

\begin{Thm} \label{T:least}
{\rm (ZFC)} For any cardinal $\kappa\geq \aleph_0$, it is consistent to believe that there is no minimal access ideal for $2^\kappa$.
\end{Thm}

{\em Proof.} Given $\kappa\geq \aleph_0$, it is well-known that the following is consistent  with ${\rm ZFC}$: $\kappa<{\rm cf}(\lambda)<\lambda$ and $2^\kappa=\lambda$. Fixing  any singular cardinal $\lambda$ with $\kappa<{\rm cf}(\lambda)<\lambda$, if $\mbox{\goth M}$ is a model of ${\rm ZFC}$ obtained by a $\kappa$-cc forcing which makes $2^\kappa=\lambda$, then   $\lambda$ remains  singular in $\mbox{\goth M}$. So it follows from Theorem \ref{T:nominimal} that in $\mbox{\goth M}$ no minimal access ideal exists for $2^\kappa$.   \telos

\section{Algebraic similarity of ${\cal P}(\kappa)$ and ${\cal C}(V)$   equipped with  ideals}
Given infinite cardinals $\lambda\leq \kappa$, let $\langle {\cal P}(\kappa),[\kappa]^{<\lambda}\rangle$ and $\langle {\cal C}(V),V\rangle$ be the Boolean algebras of all  subcollections of $\kappa$ and $V$, equipped  with the ideals $[\kappa]^{<\lambda}$ and $V$, respectively, representing some ``small sets/classes'' of the algebras.  The language of the structures is  $L_b(S)=L_b\cup\{S\}$, where  $L_b$ is the language of Boolean algebras, and  $S(\cdot)$ is a unary predicate symbol for the elements of the ideals. The purpose of this section is to prove the following.
\begin{Thm} \label{T:ideal}
{\rm (GBC)} For all infinite cardinals $\lambda\leq \kappa$,
$\langle {\cal P}(\kappa), [\kappa]^{<\lambda}\rangle\equiv_{L_b(S)}$
$\langle {\cal C}(V),V\rangle$. In particular for  $\kappa=2^\lambda$ we get $\langle{\cal P}(2^\lambda), [2^\lambda]^{\leq \lambda}\rangle \equiv_{L_b(S)}\langle {\cal C}(V),V\rangle$.
\end{Thm}

It is already known that any two powerset Boolean algebras  are elementarily equivalent. This is an easy corollary of a general result concerning algebras with the same ``invariants'' (see \cite[p. 296]{Ko89}, \S 18). However the algebras here are endowed with ideals so a stronger argument is needed. The method of proof is an application of the Ehrenfeucht-Fra\"{\i}ss\'{e} games (see e.g.  \cite{Ma02}, pp. 52f). We recall briefly some definitions.

Let $L$ be a first-order language and  let ${\cal M}=\langle M,\ldots\rangle$, ${\cal N}=\langle N,\ldots\rangle$ be two $L$-structures. We define first what an infinite game is  between ${\cal M}$, ${\cal N}$. This is a two-player game denoted $G_\omega({\cal M},{\cal N})$. The players, called I and II, play at stages $1,2,\ldots$, by choosing at each stage $i$, alternately elements $a_i\in M$ and $b_i\in N$. A play consists of $\omega$ stages. Player I may choose an element either from $M$ or from $N$, and then II must respond by choosing an element from the other structure. Thus if I chooses  $a_i\in M$ (resp. $b_i\in N$), then II must choose  $b_i\in N$ (resp. $a_i\in M$). Player II wins the play if $f=\{\langle a_i,b_i\rangle:i=1,2,\ldots\}$ is  the graph of a partial embedding from ${\cal M}$ into ${\cal N}$. Otherwise I wins. A strategy for player II in $G_\omega({\cal M},{\cal N})$ is a function $\tau$ such that for each $i$, if $u_1,\ldots,u_i$ are the first $i$ moves of player I, then the $i$-th move of player II is $\tau(u_1,\ldots,u_i)$. The strategy $\tau$ is a {\em winning} strategy for II, if for any sequence $u_1,u_2,\ldots$ of moves of player I, II wins the play if he/she follows the strategy $\tau$.

Now an Ehrenfeucht-Fra\"{\i}ss\'{e} (E-F) game between ${\cal M}$, ${\cal N}$ is a  game of the above general form, except that it  consists only of {\em finitely} many stages instead of infinite ones. It is  denoted $G_n({\cal M}, {\cal N})$ if it stops after the $n$-th stage. The winning strategy for player II or I  is defined as before. E-F games can be used to characterize elementary equivalence of two structures. Specifically the following holds.

\begin{Thm} \label{T:marker}
{\rm (\cite[Thm 2.4.6]{Ma02})} Let $L$ be a finite language without function symbols and let ${\cal M}$, ${\cal N}$ be $L$-structures. Then ${\cal M}\equiv {\cal N}$ if and only if player II has a winning strategy in the game $G_n({\cal M},{\cal N})$, for all $n>0$.
\end{Thm}
We shall use \ref{T:marker} to prove Theorem \ref{T:ideal}. In the proof we shall need the following simple fact which is provable in GB.

\begin{Lem} \label{L:classplit}
{\rm (GB)} Every proper class can be partitioned into two proper disjoint subclasses.
\end{Lem}

{\em Proof.} Let $X$ be a proper class and, as in Lemma \ref{L:around}, consider the class $R=\{rank(x):x\in X\}$ which is proper and can be written $R=\{\alpha_\xi:\xi\in On\}$ in its natural well-ordering. Let $R_0=\{\alpha_\xi: \xi \ \mbox{is a limit ordinal}\}$, and $R_1=\{\alpha_\xi: \xi \ \mbox{is a successor ordinal}\}$. If for $i=0,1$ we set $X_i=\{x\in X:rank(x)\in R_i\}$, $X_0$, $X_1$ is the required partition. \telos

\vskip 0.2in

{\em Proof of Theorem \ref{T:ideal}}. Fix infinite cardinals $\lambda\leq \kappa$ and let ${\cal M}=\langle {\cal P}(\kappa), [\kappa]^{<\lambda}\rangle$ and ${\cal N}=\langle {\cal C}(V),V\rangle$. (Remember that we can always think we are working in a model $\mbox{\goth K}=\langle K,{\cal K}\rangle$ of ${\rm GBC}$, with $V^{\mbox{\goth K}}=K$  and  ${\cal C}(V)^{\mbox{\goth K}}={\cal K}$, so in this case  ${\cal N}=\langle{\cal K},K\rangle$.) In order to comply with the condition  of Theorem \ref{T:marker} that the language of the structures must be  relational, we shall take the language $L_b$ of Boolean algebras to be $L_b=\{\sqsubseteq\}$ (instead of the more often used algebraic language $\{+,\cdot,-,0,1\}$). According to Theorem \ref{T:marker}, we have  to show that  for each $n$, Player II has a winning strategy in  $G_n({\cal M},{\cal N})$. In fact we can show that player II has a winning strategy not only in the finite games $G_n({\cal M},{\cal N})$, but also in the infinite game $G_\omega({\cal M},{\cal N})$. To see this it suffices to define, for every $n$, what the move of player II will be after the completion of the first $n$ stages.  Now since our ground structures ${\cal M}$, ${\cal N}$ are Boolean algebras, every finite subset $A$ of each  of the algebras generates a finite Boolean subalgebra {\em fully determined} by the set of its atoms which are the intersections of the elements of $A$. And  whenever a player chooses a new element from the ground algebra, this element expands the previous subalgebra to a larger finite one by refining the set of its atoms. Thus in our case each stage of the game is not obtained from  the previous one by the attachment of a single new pair of elements from each algebra,  but of several new pairs that correspond to the generated new atoms. So we can consider the rounds of the game as partial embeddings  $f_i:C_i\rightarrow D_i$, $i\geq 1$, where $C_i$, $D_i$ are  finite Boolean subalgebras of ${\cal P}(\kappa)$ and ${\cal C}(V)$ respectively, and for every $i$, $f_i\subsetneq f_{i+1}$. Moreover, $f_i$ respect the predicate $S(\cdot)$ by preserving the ideals, i.e., satisfying the condition

\vskip 0.1in

(*) \quad  For every $a\in C_i$, $a\in [\kappa]^{<\lambda} \Leftrightarrow f_i(a)\in V$.

\vskip 0.1in

\noindent Assume that the first $n$ rounds of the game have been completed, i.e., embeddings $f_i:C_i\rightarrow D_i$ which satisfy (*) have been defined for all $i\leq n$,  and suppose player I starts the $(n+1)$-th round  by choosing either some $a\in {\cal P}(\kappa)$ or some $b\in {\cal C}(V)$. We have to define the response of player II, i.e. a $b\in {\cal C}(V)$,  if player I chooses $a\in {\cal P}(\kappa)$, or an $a\in {\cal P}(\kappa)$,  if player I chooses $b\in {\cal C}(V)$, and the extension  $f_{n+1}:C_{n+1}\rightarrow D_{n+1}$ of $f_n$ so that  $\langle a,b\rangle\in f_{n+1}$ and $f_{n+1}$ satisfies (*). Without loss of generality we  assume that $a\in {\cal P}(\kappa)$. We shall describe the definition of $b$ and $f_{n+1}$, and this will be the description of the winning strategy of player II. A quite similar strategy would be described if player I had started with a  $b\in {\cal C}(V)$ and player II should choose $a\in {\cal P}(\kappa)$.

{\em The strategy of player II}: Let $c_1,\ldots, c_k$ be the atoms of the algebra $C_n$ which form a partition of $\kappa$. Then $f_n(c_1),\ldots, f_n(c_k)$ are the atoms of  $D_n$ which  form a partition of $V$. For  each  $i=1,\ldots,k$, let $x_i=c_i\cap a$ and $x^\prime_i=c_i\cap (\kappa\backslash a)$ and  let
$$X=\{x_i, x_i^\prime\neq \emptyset:i=1,\ldots,k\}.$$
$X$ is a new partition of $\kappa$ which refines  the partition of $c_i$'s and determines a finite subalgebra $C_{n+1}$ of ${\cal P}(\kappa)$ such that $a\in C_{n+1}\supseteq C_n$. In order to define $f_{n+1}$ and $f_{n+1}(a)=b$, clearly it suffices to define $f_{n+1}(x_i)=y_i$ and $f_{n+1}(x_i')=y_i'$ for all $x_i,x_i'\in X$, so that $y_i,y_i'$ form a suitable partition of $f_n(c_i)$. In particular we must have
\begin{equation} \label{E:extra}
|x_i|<\lambda \ (|x^\prime_i|<\lambda)\Leftrightarrow y_i\in V \ (y^\prime_i\in V).
\end{equation}
If for some $i$, $c_i\subseteq a$, or $c_i\subseteq \kappa\backslash a$, then $x_i=c_i$ or $x_i^\prime=c_i$, respectively, for which $f_n$ is already defined. So  if  $x_i=c_i$ or $x_i'=c_i$, we just let $f_{n+1}(x_i)=f_n(c_i)$ and $f_{n+1}(x_i^\prime)=f_n(c_i)$. In all other cases $x_i,x_i^\prime$ form a proper partition of $c_i$ and we  consider the following cases with respect to the cardinalities of $x_i$ and $x_i^\prime$.

1) $|x_i|, |x^\prime_i|\geq \lambda$. Then $|c_i|\geq \lambda$ so, by (*) for $f_n$, $f_n(c_i)$ is a proper class. In this case it suffices to take, using  Lemma  \ref{L:classplit}, $y_i$, $y^\prime_i$ to be  complementary proper subclasses of $f_n(c_i)$.

2) $|x_i|<\lambda$ and $|x^\prime_i|\geq\lambda$ (or vice versa). Then $|c_i|\geq\lambda$, hence  $f_n(c_i)$ is a proper class.  Choose $\emptyset\neq y_i\subseteq f_n(c_i)$ to be a set and let $y^\prime_i=f_n(c_i)\backslash y_i$. In particular, if $x_i$ is finite, we pick $y_i$ so that $|y_i|=|x_i|$.

3) $|x_i|, |x^\prime_i|<\lambda$. Then $|c_i|<\lambda$. Therefore  $f_n(c_i)$ is a set.  If $x_i,x_i'$ are both infinite, we pick any infinite  $y_i$, $y^\prime_i$ that form a partition of $f_n(c_i)$. If either $x_i$ or/and $x_i'$  is finite we take the partition $y_i,y_i^\prime$ so that  $|y_i|=|x_i|$ or/and $|y_i'|=|x_i'|$, respectively. (This is always possible, even when $c_i$ is finite. Because the same strategy has been followed by player II at the previous stages, so if $c_i$ is finite  the  strategy has yielded  $f_n(c_i)$ such that $|f_n(c_i)|=|c_i|$.) Let
$$Y=\{y_i,y^\prime_i:i=1,\ldots,k\}.$$
$Y$ is the set of atoms of a  subalgebra $D_{n+1}\supseteq D_n$ of ${\cal C}(V)$ and the function  $f_{n+1}:C_{n+1}\rightarrow D_{n+1}$ for which  $f_{n+1}(x_i)=y_i$ and $f_{n+1}(x_i^\prime)=y_i^\prime$ induces  a partial isomorphism  that extends $f_n$. Thus clearly   $f_{n+1}$ satisfies (\ref{E:extra}). It satisfies also  condition (*)  above. Because by (\ref{E:extra}),  $f_{n+1}$ already satisfies (*) for the   atoms of $C_{n+1}$. If now  $x$ is any element of $C_{n+1}$  and $z_1,\ldots, z_m$ are the atoms of $C_{n+1}$ below $x$, then $x=z_1\cup\cdots\cup z_m$, hence $f_{n+1}(x)=f_{n+1}(z_1)\cup\cdots\cup f(z_m)$. So clearly  $|x|<\lambda$ iff $|z_i|<\lambda$ for all $i=1,\ldots,m$, iff $f_{n+1}(z_i)$ is a set for all $i=1,\ldots,m$, iff $f_{n+1}(x)$ is a set. Thus $f_{n+1}$ satisfies (*).

This completes the description of the strategy  of player II. It remains to prove that the strategy is winning, that is, all these  $f_i$ are  partial embeddings from ${\cal M}$ into ${\cal N}$. Namely it suffices to show  the following:

\vskip 0.1in

{\em Claim}. Let $F$ be the set of all finite functions  from ${\cal M}$ into ${\cal N}$ which are defined through  E-F games, in which player II follows the above described strategy. Then every $f\in F$ is a partial embedding, i.e., for every  $\phi(x_1,\ldots,x_n)$ of $L_b(S)$ and every  $a_1,\ldots,a_n\in {\cal P}(\kappa)$ such that  $a_1,\ldots,a_n\in dom(f)$,
\begin{equation} \label{E:entails}
{\cal M} \models\phi(a_1,\ldots,a_n) \Leftrightarrow {\cal N}\models \phi(f(a_1),\ldots,f(a_n)).
\end{equation}

{\em Proof of Claim.} We prove that (\ref{E:entails}) holds for every $f\in F$ and for every $\phi$   by induction on the length of $\phi$. The atomic formulas of $L_b(S)$ are $x\sqsubseteq y$ and $S(x)$ and clearly the ordering $\sqsubseteq$ is preserved by all $f\in F$. Also, we saw above that every $f\in F$  satisfies condition (*), so for every $a\in dom(f)$
$${\cal M}\models S(a) \Leftrightarrow a\in [\kappa]^{<\lambda}\Leftrightarrow  f(a)\in V\:\Leftrightarrow {\cal N}\models S(f(a)).$$
The induction  steps for the connectives are straightforward for all $f\in F$, and the only nontrivial step of the induction is the  quantifier step. So given a $\phi(x,x_1,\ldots,x_n)$ we assume that all $f\in F$ satisfy  (\ref{E:entails}) for $\phi$, and prove this  for $\exists x\phi(x,x_1,\ldots,x_n)$, i.e., for all $a_1,\ldots,a_n$ in $dom(f)$
\begin{equation} \label{E:step}
{\cal M}\models\exists x\phi(x,a_1, \ldots,a_n)\Leftrightarrow {\cal N}\models\exists x\phi(x,f(a_1),\ldots,f(a_n)).
\end{equation}
Suppose  ${\cal M}\models\exists x\phi(x,a_1, \ldots,a_n)$ for given $a_1,\ldots,a_n$ in $dom(f)$. Then ${\cal M}\models \phi(a,a_1, \ldots,a_n)$ for some $a\in {\cal P}(\kappa)$. Now, by definition,  $f$  is  the $n$-th stage of some E-F game between ${\cal M}$ and ${\cal N}$, so $a$ can be the move of player I at the $(n+1)$-th stage of the game. Then player II responds following the above strategy by picking a $b\in {\cal C}(V)$ and extending $f$ to $g\supseteq f\cup\{\langle a,b\rangle\}$ so that $g\in F$. By our  induction assumption $g$ satisfies (\ref{E:entails}) for $\phi$, hence  ${\cal N}\models\phi(g(a),g(a_1),\ldots,g(a_n))$, or ${\cal N}\models\phi(b,f(a_1),\ldots,f(a_n))$, therefore ${\cal N}\models\exists x\phi(x,f(a_1),\ldots,f(a_n))$. This proves ``$\Rightarrow$'' of (\ref{E:step}). For the converse assume ${\cal N}\models \exists x\phi(x,f(a_1),\ldots,f(a_n))$ and ${\cal N}\models\phi(b,f(a_1),\ldots,f(a_n))$ for a $b\in {\cal C}(V)$. As before $b$ can also be the move of player I at the $(n+1)$-th stage, and then player II responds by choosing an $a\in {\cal P}(\kappa)$ and $g\supseteq f\cup\{\langle a,b\rangle\}$ so that $g$ belongs to $F$. Then ${\cal N}\models\phi(g(a),g(a_1),\ldots,g(a_n))$ and by the induction assumption ${\cal M}\models \phi(a,a_1,\ldots,a_n)$, whence  ${\cal M}\models \exists x\phi(a_1,\ldots,a_n)$. This proves the Claim and completes the proof of the Theorem. \telos

\vskip 0.2in

The similarity of structures $\langle {\cal P}(\kappa), [\kappa]^{<\lambda}\rangle$ and $\langle {\cal C}(V),V\rangle$ expressed  by Theorem \ref{T:ideal} is of a rather algebraic character, since it is about Boolean algebras equipped with ideals. However it can be transcribed into one with more set-theoretic flavor, i.e., one using the language of set theory. Specifically, the elements of the ground set $\kappa$ and class $V$ can be treated as  {\em atoms} (or {\em urelements}), while ${\cal P}(\kappa)$ and ${\cal C}(V)$ as collections of sets/classes  of atoms.   So  let $L_1(S)$ be the language consisting of $\in$, $S(\cdot)$, variables $p,q,\ldots$ for urelements, and variables  $x,y,\ldots$ for sets. The atomic formulas of $L_1(S)$ are: $p=q$, $p\in x$, and $S(x)$.  The intended $L_1(S)$-structures are $\langle\kappa,{\cal P}(\kappa),[\kappa]^{<\lambda}\rangle$, for sets, and  $\langle V,{\cal C}(V),V\rangle$ for the universe. However a  technical problem here is the fact that the  two kinds of entities, atoms and sets, should  be disjoint, while  $\kappa\cap {\cal P}(\kappa)\neq \emptyset$ (and similarly for $V$). We can bypass this obstacle by working with suitable copies of $\kappa$ and $V$. Namely we set
$$\kappa^*=\{\langle x,0\rangle:x\in\kappa\} \ \mbox{and} \ V^*=\{\langle x,0\rangle:x\in V\}.$$
Then $\kappa^*\cap {\cal P}(\kappa^*)=\emptyset$ and $V^*\cap {\cal C}(V^*)=\emptyset$. Under this coding the $L_b(S)$ structures $\langle{\cal P}(\kappa),[\kappa]^{<\lambda}\rangle$ and  $\langle{\cal C}(V),V\rangle$ become isomorphic  to $\langle{\cal P}(\kappa^*),[\kappa^*]^{<\lambda}\rangle$ and $\langle{\cal C}(V^*),{\cal P}(V^*)\rangle$, respectively, where ${\cal P}(V^*)$ is the class of subsets of $V^*$. Then in view of Theorem \ref{T:ideal} we have: For all infinite cardinals $\lambda\leq\kappa$
\begin{equation} \label{E:copies}
\langle{\cal P}(\kappa^*),[\kappa^*]^{<\lambda}\rangle\equiv_{L_b(S)} \langle{\cal C}(V^*),{\cal P}(V^*)\rangle.
\end{equation}
Based on (\ref{E:copies}) we can easily deduce the following.

\begin{Thm} \label{T:reduct}
For all  infinite $\lambda\leq \kappa$,
$$\langle\kappa^*,{\cal P}(\kappa^*),[\kappa^*]^{<\lambda}\rangle\equiv_{L_1(S)}  \langle V^*,{\cal C}(V^*),{\cal P}(V^*)\rangle.$$
\end{Thm}

{\em Proof}. Define an interpretation $\phi\mapsto \phi^+$ from the atomic formulas of $L_1(S)$ to formulas of $L_b(S)$ as follows: Fix enumerations  $p_n, x_n$, $n\in \N$, of the variables of $L_1(S)$ for atoms and sets, respectively, and  an enumeration $y_n$, $n\in \N$,  of the variables of $L_b(S)$. The atomic formulas of $L_1(S)$ are $p_i=p_j$, $p_i\in x_j$ and $S(x_j)$. Define

(a) $(p_i=p_j)^+:=(y_{2i}=y_{2j})$,

(b) $(p_i\in x_j)^+:=Atom(y_{2i}) \wedge y_{2i}\sqsubseteq y_{2j+1}$,

(c) $S(x_j)^+=S(y_{2j+1})$, \\
where $Atom(y)$ is the formula of $L_b(S)$ saying that ``$y$ is an atom''. This interpretation is also extended to all formulas of $L_1(S)$ as expected. (As follows from clauses (a)-(c) above, the variable $p_i$ translates to $y_{2i}$ while $x_j$ translates to $y_{2j+1}$. This translation determines also of course the interpretation of quantifiers.)  Then it is easy to see that for all  $\phi\in L_1(S)$,
$$\langle\kappa^*,{\cal P}(\kappa^*),[\kappa^*]^{<\lambda}\rangle\models \phi \ \Leftrightarrow \ \langle{\cal P}(\kappa^*),[\kappa^*]^{<\lambda}\rangle\models \phi^+,$$
and
$$\langle V^*,{\cal C}(V^*),{\cal P}(V^*)\rangle\models \phi \Leftrightarrow \langle{\cal C}(V^*),{\cal P}(V^*)\rangle\models\phi^+.$$
By (\ref{E:copies}),
$$\langle{\cal P}(\kappa^*),[\kappa^*]^{<\lambda}\rangle\models \phi^+\Leftrightarrow
\langle{\cal C}(V^*),{\cal P}(V^*)\rangle\models\phi^+,$$
therefore
$$\langle\kappa^*,{\cal P}(\kappa^*),[\kappa^*]^{<\lambda}\rangle\models \phi \ \Leftrightarrow \langle V^*,{\cal C}(V^*),{\cal P}(V^*)\rangle\models \phi.$$
\telos

\vskip 0.2in

In fact Theorems \ref{T:ideal} and \ref{T:reduct} are essentially equivalent, since a proof of \ref{T:ideal} from \ref{T:reduct} can also be obtained by defining a converse  interpretation $\phi\mapsto \phi^\prime$ from atomic formulas of $L_b(S)$ to formulas of $L_1(S)$ as follows:

$(x\sqsubseteq y)^\prime:=\forall p (p\in x\rightarrow p\in y)$,

$(x=y)^\prime:=(x=y)$, and

$S(x)^\prime:=S(x)$. \\
This interpretation, extended to all formulas, yields: For every $\phi$ of $L_b(S)$,
$$\langle{\cal P}(\kappa^*),[\kappa^*]^{<\lambda}\rangle\models \phi \  \Leftrightarrow \ \langle\kappa^*,{\cal P}(\kappa^*),[\kappa^*]^{<\lambda}\rangle\models \phi^\prime$$ and
$$\langle{\cal C}(V^*),{\cal P}(V^*)\rangle\models \phi \Leftrightarrow \langle V^*,{\cal C}(V^*),{\cal P}(V^*)\rangle\models \phi^\prime.$$
Thus  \ref{T:reduct} implies \ref{T:ideal}.

\section{Order-theoretic similarity of uncountable cardinals with  the class of ordinals}
By the axiom of choice every uncountable set $A$ can be identified with an uncountable cardinal $\kappa$,  and by the axiom of global choice $vN$ the universe $V$ can be identified with   the class of ordinals $On$. However  these identifications  {\em add structure} to $A$ and $V$, namely the well-orderings of $\kappa$ and $On$, which we denote by the same symbol $<$. Admitting the naturalness of choice principles, it is then reasonable to compare  the first-order structures $\langle \kappa,<\rangle$, $\langle \lambda,<\rangle$  for uncountable $\kappa,\lambda$, as well as  $\langle \kappa,<\rangle$ and $\langle On,<\rangle$.  In this section we shall show how the main result of \cite{Tz04a} (Theorem 2.1), which says that the monadic $\forall_1^1$ positive theory of linear ordering is the same for all uncountable cardinals,  can be extended  also  over $\langle On,<\rangle$.

The starting point  of \cite{Tz04a} was the following result, which is a special case of  \cite[Cor. 44]{DMT78}:

\begin{Thm} \label{T:arxhpalia}
{\rm (\cite[Thm 1.1]{Tz04a})} For all uncountable cardinals $\kappa$, $\lambda$, 
$\langle \kappa,<\rangle\equiv_{L_{ord}}\langle\lambda,<\rangle$, where $L_{ord}=\{\prec\}$ is the first-order language of ordering.
\end{Thm}
Remember that the ordering $<$ of ordinals is a well-ordering, that is a second-order property, while the equivalence $\equiv_{L_{ord}}$  refers only to the first-order properties of $<$. Theorem \ref{T:arxhpalia}  is in fact a consequence of Corollary 44 of  \cite{DMT78} and, essentially, the same argument that we  used in order to justify it in \cite{Tz04a} can be applied also to justify  the following extension.

\begin{Thm} \label{T:arxhnea}
For every uncountable cardinal $\kappa$,
$\langle \kappa,<\rangle\equiv_{L_{ord}}\langle On,<\rangle$.
\end{Thm}

{\em Proof.} (Sketch)  By Corollary  44 of  \cite{DMT78}, for any ordinals $\alpha$, $\beta$, $\langle \alpha,<\rangle\equiv \langle \beta,<\rangle$ if and only if $\alpha$, $\beta$ are {\em congruent modulo $\omega^\omega$}, i.e.,  there are $\xi$, $\eta$, $\delta$ such that  $\delta<\omega^\omega$,  $\alpha=\omega^\omega\cdot\xi+\delta$, $\beta=\omega^\omega\cdot\eta+\delta$ and either $\xi=\eta=0$ or $\xi\neq 0$ and $\eta\neq 0$.\footnote{Important notice: when I was preparing \cite{Tz04a}, Professor Doner had kindly informed me in a private communication,  that the  definition of {\em congruence modulo $\omega^\omega$} as given in \cite{DMT78}, p. 51, is mistaken. The correct definition is that given above and can be found in \cite[p. 6]{Do72}.}   Since for all cardinals $\kappa,\lambda>\omega$, $\kappa=\omega^\omega\cdot\kappa$ and $\lambda= \omega^\omega\cdot\lambda$, it follows that $\kappa, \lambda$ are congruent modulo $\omega^\omega$. This already proves Theorem \ref{T:arxhpalia}.

Concerning the justification of the present claim, without any risk of fallacy we can treat  $On$ in this particular situation as the  ``greatest (limit) ordinal'',  and extend the usual (right) multiplication of ordinals over $On\cup\{On\}$ in the obvious way. Specifically, just as for any $\alpha$ and any limit $\beta$ we have  $\alpha\cdot\beta=\sup\{\alpha\cdot\gamma:\gamma<\beta\}$, simply define  $\alpha\cdot On=\sup\{\alpha\cdot\gamma:\gamma<On\}$. Then trivially $\omega^\omega\cdot On=On$,  so  for every uncountable cardinal $\kappa$, $\kappa$ and $On$ are congruent modulo $\omega^\omega$. Applying again Corollary  44 of  \cite{DMT78} to $\langle \kappa,<\rangle$ and $\langle On,<\rangle$, we are done.  (As for the details  why the relation of  congruence modulo $\omega^\omega$ entails the elementary equivalence of well-ordered sets or classes,  one should consult  \cite{DMT78}, especially  the theory ${\cal W}$ (Def. 2) of ``weak well-orderings'', the notion of ``canonical order type'' (Def. 39),  and Theorem 41 which relates these notions with elementarily equivalent ordinals.)  \telos

\vskip 0.2in

We proceed next to extend the main theorem of \cite{Tz04a}, which establishes a similarity significantly stronger  than that of Theorem \ref{T:arxhpalia}, so as to hold between every  uncountable cardinal and $On$. We need first some terminology and notation.

We  first extend the language of ordering $L_{ord}=\{\prec\}$ used previously for the structures $\langle \kappa,<\rangle$ (with lower-case individual variables $x,y,z,\ldots$) to the {\em monadic} (second-order) language $L_{mon}$ which is $L_{ord}$ augmented with $\in$ and upper-case set variables $X,Y,Z,\ldots$.  A formula of $L_{mod}$ without quantified set variables is called {\em normal}. A formula of the form $\forall X\phi$   (resp. $\exists X\phi$), where $\phi$ is normal is called a $\forall_1^1$ (resp. $\exists_1^1$) formula. The standard interpretations of ${\cal L}_{mon}$ are the structures $\langle A,<, {\cal P}(A)\rangle$ where $\langle A,<\rangle$ is an ordered set that interprets the formulas of $L_{ord}$,  and ${\cal P}(A)$ is the range of set variables.  However there are more general interpretations. They are of the form $\langle A,<,{\cal A}\rangle$, where  ${\cal A}\subset {\cal P}(A)$. ${\cal A}$  may be some class of ``small'' sets, i.e.,  an ideal of ${\cal P}(A)$, or, in the opposite direction, a class of ``large'' sets. Most often  these large sets are  the sets cofinal in $\langle A,<\rangle$, so let us set
$${\rm Cof}(A)=\{X\subseteq A: \mbox{$X$ is cofinal in
$\langle A,<\rangle$}\}.$$
$A$ can be also a proper class and ${\cal A}$ be a class of subclasses of $A$.
When we deal with an  ordered collection  which is a  proper class, like $On$, then, in order to interpret the set variables $X$, $Y$, we must work again  in a class theory.  Then as usual we write  ${\cal C}(On)$   for the collection of all subclasses of $On$. In particular below we shall deal with the analogue of ${\rm Cof}(\kappa)$ above, which is
$${\rm Cof}(On)=\{X\subseteq On:\mbox{$X$ is cofinal in} \
\langle On,<\rangle\}$$
and consists of the proper subclasses of $On$.

\begin{Def} \label{D:positive}
{\em Let $\phi(X)$ be a normal formula of ${\cal L}_{mon}$ with at most one (free) set variable $X$.  $\phi(X)$ is said to be} positive in $X$ {\em  (or just} positive {\em) if it belongs to the smallest class of formulas which (a) contains all formulas of $L_{ord}$, (b) contains the atomic formulas $x\in X$ and (c) is closed under the logical operations $\wedge$,  $\vee$, and the first-order quantifiers $\exists$ and $\forall$. }
\end{Def}

The main result  of \cite{Tz04a} is  that the structures $\langle \kappa,<,{\rm Cof}(\kappa)\rangle$, for all uncountable cardinals $\kappa$, are elementarily equivalent with respect to $\forall^1_1$ positive formulas of $L_{mod}$.

\begin{Thm} \label{T:mainpalio}
{\rm (\cite[Thm. 2.1]{Tz04a})} In ${\rm ZFC}$ the following holds: For all uncountable cardinals $\kappa,\lambda$,
$$\langle \kappa,<,{\rm Cof}(\kappa)\rangle\equiv_{pos}^{\forall_1^1} \langle \lambda,<,{\rm Cof}(\lambda)\rangle,$$
i.e., for every positive $\phi(X)$,
$$\langle \kappa,<,{\rm Cof}(\kappa)\rangle\models(\forall X)\phi(X)\Leftrightarrow \langle \lambda,<,{\rm Cof}(\lambda)\rangle\models (\forall X)\phi(X).$$

\end{Thm}

Here we strengthen \ref{T:mainpalio} to the following.

\begin{Thm} \label{T:maineo}
{\rm (${\rm GBC}$)}  For every  uncountable cardinal $\kappa$,
$$\langle \kappa,<,{\rm Cof}(\kappa)\rangle\equiv_{pos}^{\forall_1^1} \langle On,<,{\rm Cof}(On)\rangle.$$
\end{Thm}

{\em Proof.} (Sketch)   The proof goes exactly along the lines of the proof of \ref{T:mainpalio} and is based on the  {\em syntactic form} of  positive (normal) formulas of $L_{mod}$ which has been discovered  by Y. Moschovakis. Namely, every formula $\phi(X)$ positive in $X$ is equivalent  over any structure $\langle A,<,{\cal P}(A)\rangle$,  for all  $X\neq A$,  to a formula having  prenex form  $(\overline{Q}\overline{w})(\forall u)(\theta(\overline{w},u)  \rightarrow  u\in X)$, where $\overline{w}=\langle w_1,\ldots,w_n\rangle$ is a string of variables, $\overline{Q}=\langle Q_1,\ldots,Q_n\rangle$ a string of quantifiers and $\theta(\overline{w},u)$ is  quantifier-free. (Note that, since $\phi(X)$ does not contain set quantifiers,  the above prenex  form for $\phi(X)$ is valid also over every  structure $\langle A,<,{\cal A}\rangle$, with  ${\cal A}\subseteq {\cal P}(A)$.) In view of the preceding  syntactic characterization of $\forall_1^1$ positive formulas, the proof splits  into several cases  that correspond to the various  forms of the string  $\overline{Q}\overline{w}$, as well as to  the forms  of $\theta$ as a Boolean combination of atoms and negated atoms. Then, through an exhaustive inspection of all these  cases, it is shown that we can essentially  reduce the truth of $(\forall X)\phi(X)$ to  the truth of a {\em first-order} formula. But this  means, in view of Theorem \ref{T:arxhnea},  that $(\forall X)\phi(X)$ is absolute between $\langle \kappa,<,{\rm Cof}(\kappa)\rangle$ and $\langle On,<,{\rm Cof}(On)\rangle$. \telos

\vskip 0.2in

If we work in ${\rm GB}+vN$, where a bijection  between $On$ and $V$ exists, the well-ordering of $On$ is transferred to a well-ordering of $V$ so the following holds.

\begin{Cor} \label{C:universe}
In ${\rm GB}+vN$ it is provable that there is a well-ordering $\prec$ of $V$ such that for every uncountable cardinal $\kappa$,
$$\langle \kappa,<,{\rm Cof}(\kappa)\rangle\equiv_{pos}^{\forall_1^1} \langle V,\prec,{\rm Cof}(V)\rangle.$$
\end{Cor}

\vskip 0.1in

Concerning the optimality of Theorem \ref{T:maineo} with  respect to the indistinguishability of structures  $\langle \kappa,<,{\rm Cof}(\kappa)\rangle$ and $\langle On,<,{\rm Cof}(On)\rangle$, and the analogous question about structures $\langle \kappa,<,{\rm Cof}(\kappa)\rangle$ and $\langle \lambda,<,{\rm Cof}(\lambda)\rangle$ for uncountable $\kappa$, $\lambda$, the reader is referred to the corresponding discussion in \cite{Tz04a}: It is not known  whether  $\langle \kappa, <,{\rm Cof}(\kappa)\rangle $ and $\langle \lambda,<, {\rm Cof}(\lambda)\rangle $ are indistinguishable with respect to all $\forall^1_1$ and $\exists^1_1$ formulas.  What we know is that the structures $\langle \omega_1,<, {\rm Cof}(\omega_1)\rangle$ and $\langle \omega_2, <,{\rm Cof}(\omega_2)\rangle$ can be  distinguished by a sentence  $\phi_{\omega_1}$ of the form $(\forall x)(\exists X)\psi$, where $\psi$ is normal but not positive (see section 5 of \cite{Tz04a}). This is easily generalized to all $\omega_m$, $\omega_n$. A fortiori $\phi_{\omega_1}$ distinguishes also the structures   $\langle \omega_1,<, {\rm Cof}(\omega_1)\rangle$ and $\langle On,<,{\rm Cof}(On)\rangle$.

\vskip 0.2in

\textbf{Acknowledgement} My thanks go to an anonymous referee whose remarks helped me fix some unclarities and improve considerably the presentation of the article.


\begin{thebibliography}{99}
\bibitem{Do72}
   J.E. Doner, Definability in  extended arithmetic of ordinal numbers,  {\em
   Dissertationes Math.  (Rozprawy Mat.)} {\bf 96} (1972),   50 pp.
\bibitem{DMT78}
   J.E. Doner, A. Mostowski and A. Tarski, The elementary theory of well-ordering --
   a metamathematical study, in: {\em Logic Colloquium '77,} A. Macintyre, L.
   Pacholski and  J. Paris (eds.), North Holland, 1978, pp. 1-54.
\bibitem{Ko89}
   S. Koppelberg, Metamathematics (of Boolean a lgebras), in: {\em Handbook of
   Boolean Algebras},  J.D. Monk and R. Bonnet (eds.), vol. 1, Elsevier Science
   Publishers B.V. 1989, pp. 287-307.
\bibitem{Ma02}
   D. Marker, {\em Model Theory: An Introduction}, Springer, 2002.
\bibitem{Pa71}
   R. Parikh, Existence and feasibility in arithmetic, {\em J. Symb. Logic}
  \textbf{36} (1971), no. 3, 494-508.
\bibitem{Tz03}
   A. Tzouvaras, Positive set-operators of low complexity,  {\em Math. Log. Quar.}
   {\bf 49} (2003),  no. 3,  284-292.
\bibitem{Tz04a}
   A. Tzouvaras, Uncountable cardinals have the same monadic $\forall_1^1$-positive
   theory over large sets, {\em Fund. Math.} \textbf{181} (2004), no. 2, 125-142.
\bibitem{Tz04b} A. Tzouvaras, What is so special with the power set operation?,  {\em
   Arch.  Math. Logic}  {\bf 43} (2004),  no. 6, 723-737.
\bibitem{Tz10}
   A. Tzouvaras, Localizing the axioms, {\em Arch. Math. Logic} {\bf 49} (2010),
   no. 5, 571-601.
\end{thebibliography}
\end{document}